# MIXED POISSON APPROXIMATION OF NODE DEPTH DISTRIBUTIONS IN RANDOM BINARY SEARCH TREES

By Rudolf Grübel and Nikolče Stefanoski

*Universität Hannover*

We investigate the distribution of the depth of a node containing a specific key or, equivalently, the number of steps needed to retrieve an item stored in a randomly grown binary search tree. Using a representation in terms of mixed and compounded standard distributions, we derive approximations by Poisson and mixed Poisson distributions; these lead to asymptotic normality results. We are particularly interested in the influence of the key value on the distribution of the node depth. Methodologically our message is that the explicit representation may provide additional insight if compared to the standard approach that is based on the recursive structure of the trees. Further, in order to exhibit the influence of the key on the distributional asymptotics, a suitable choice of distance of probability distributions is important. Our results are also applicable in connection with the number of recursions needed in Hoare's [*Comm. ACM* **4** (1961) 321–322] selection algorithm Find.

## 1. Introduction.

The classical algorithm for storing data sequentially into a binary search tree proceeds as follows: The first item is put into the root node; subsequent elements are compared to the existing nodes, starting with the root, moving to the left if smaller than and to the right if greater than the content of the node until an external node is found. If there are $n$ distinct (and comparable) values, then we obtain a random binary tree if we assume that all permutations of the data are equally likely. This data structure and its properties are discussed in the standard texts of the area; see, for example, Knuth (1973), Cormen, Leiserson and Rivest (1990) and Sedgewick and Flajolet (1996). Mahmoud (1992) gives a book-length treatment of random search trees.

Suppose now that a binary search tree is associated with a random permutation of the set $\{1,2,\dots,n\}$ in the above manner. One of the quantities of interest in this structure is the depth $X_{n,l}$ of the node containing $l$, that is, its distance from the root; $1+X_{n,l}$ is the number of steps needed to retrieve the value $l$ ("successful search"). Arora and Dent (1969), in an early paper on the subject, obtained a simple and explicit formula for the corresponding expectation,

$$\text{(AD)} \qquad E(1+X_{n,l}) = H_l + H_{n+1-l} - 1,$$

where $H_k := \sum_{i=1}^{k} 1/i$, $k \in \mathbb{N}$, are the harmonic numbers. This result implies that the average number of steps needed grows logarithmically only. It is easily seen, however, that $X_{n,l}$ can be as large as $n-1$, which motivates a closer analysis of its distribution.

In contrast to many other characteristics of the tree such as its height or total path length, the depth depends on two parameters, the size $n$ of the base set and the key value (or label) $l$ of the node, which complicates the analysis. Averaging the distributions over the second parameter avoids this problem; the result can be interpreted as the distribution of the depth of a key or node selected uniformly at random from the available range $\{1,\dots,n\}$. Louchard (1987) obtained a corresponding asymptotic normality result; see also Section 2.5 in Mahmoud (1992). The distance of two randomly selected nodes has recently been investigated by Mahmoud and Neininger (2003). Averaging leads to a loss of information, though. For example, it is immediate from (AD) that

$$\lim_{n\to\infty} \frac{EX_{n,1}}{\log n} = 1, \qquad \lim_{n\to\infty} \frac{EX_{n,\lceil n/2\rceil}}{\log n} = 2,$$

that is, the depth of the node with the smallest key is only about half of that of the node with the median key value on average, if the size of the base set is large.

Our intention here is to obtain distributional approximations and asymptotics for $X_{n,l}$ that are sufficiently precise to show the dependence of the depth of a node on its key. The main tool is a distributional representation of $X_{n,l}$ in terms of mixed and compounded distributions from well-known families (Theorem 1). In contrast to many investigations in this area we do not base our analysis on a recursion for the quantities of interest, but exploit the relationship to records which seems to have been noticed first by Devroye (1988). Devroye used this connection to investigate the depth of the last node; he wrote that it "allows us to obtain …hopefully insightful proofs … ." The representation can also be used to obtain the expectation of $X_{n,l}$ and therefore leads to an alternative proof for Arora and Dent's (1969) formula. Somewhat to our surprise, asymptotic normality in the sense that $(X_{n,l_n} - EX_{n,l_n})/\sqrt{EX_{n,l_n}}$ converges in distribution to a standard normal variable holds for *every* sequence $(l_n)_{n\in\mathbb{N}}$. This result has also been obtained by Devroye and Neininger (2004). It implies Louchard's (1987) result for randomly selected nodes, but it can also be used to see the influence of the key on the node depth on the level that is apparent from the consequence of (AD) mentioned above: If the key value $l_n$ varies with $n$ such that $\log l_n/\log n \to t \in [0,1]$, then $(X_{n,l_n} - (1+t)\log n)/\sqrt{(1+t)\log n}$ is asymptotically standard normal. However, if $l_n/n \to t$ as $n \to \infty$, then the approximating normal distribution does not depend on $t$, as long as $0 < t < 1$. Hence, with this level of detail only extreme values of the key will have a noticeable influence on the depth distribution.

The proof of asymptotic normality is based on a Poisson approximation result (Theorem 3), where we use total variation distance. If we replace the total variation distance by an appropriate Wasserstein metric, then a mixed Poisson approximation is needed since with this metric shifts are not swamped by the fact that $EX_{n,l_n} \to \infty$ as $n \to \infty$. Indeed, the mixing distribution will asymptotically be close to a shifted and reflected exponential distribution, with shift $2\log n + 2\gamma + \log(t(1-t))$ depending on $t := \lim_{n\to\infty} l_n/n$ (Theorem 6; $\gamma$ denotes Euler's constant).

These results are given in the next section. In the final section we discuss various consequences of our results and also relate these to the number of recursions needed by Hoare's (1961) selection algorithm FIND.

We write $\mathcal{L}(X)$ for the distribution of the random variable $X$, with $X \stackrel{\text{distr}}{=} Y$ abbreviating $\mathcal{L}(X) = \mathcal{L}(Y)$, and $\mathbb{1}_A$ for the indicator function of the set $A$. Instead of $\mathcal{L}(X) = \mu$, with some probability distribution $\mu$, we also write $X \sim \mu$. Distributional convergence is denoted by $\stackrel{\text{distr}}{\to}$ and $N(0,1)$ is the standard normal distribution, so that $X_n \stackrel{\text{distr}}{\to} Z$, $Z \sim N(0,1)$ is short for

$$\lim_{n\to\infty} P(X_n \le x) = \Phi(x) := \frac{1}{\sqrt{2\pi}} \int_{-\infty}^{x} e^{-y^2/2}\, dy \qquad \text{for all } x \in \mathbb{R}.$$

**2. Results.** Our first result displays the distribution of $X_{n,l}$ in terms of mixed and compounded standard distributions from the Bernoulli, uniform and hypergeometric families. The representation becomes transparent once we consider an example. Suppose we have $n = 20$ and $l = 11$. A particular permutation is given in the first line of Table 1.

In the second line of Table 1 the part of the permutation to the left of the element of interest is divided into those that are greater $(+)$ or smaller $(-)$ than this element. The third line marks the descending $(\downarrow)$ and ascending $(\uparrow)$ records in these sublists, where the $i$th element $x_i$ of a list $(x_1, \ldots, x_n)$ of numbers is a descending record if $x_i = \min_{1\le j\le i} x_j$, ascending if $x_i = \max_{1\le j\le i} x_j$.

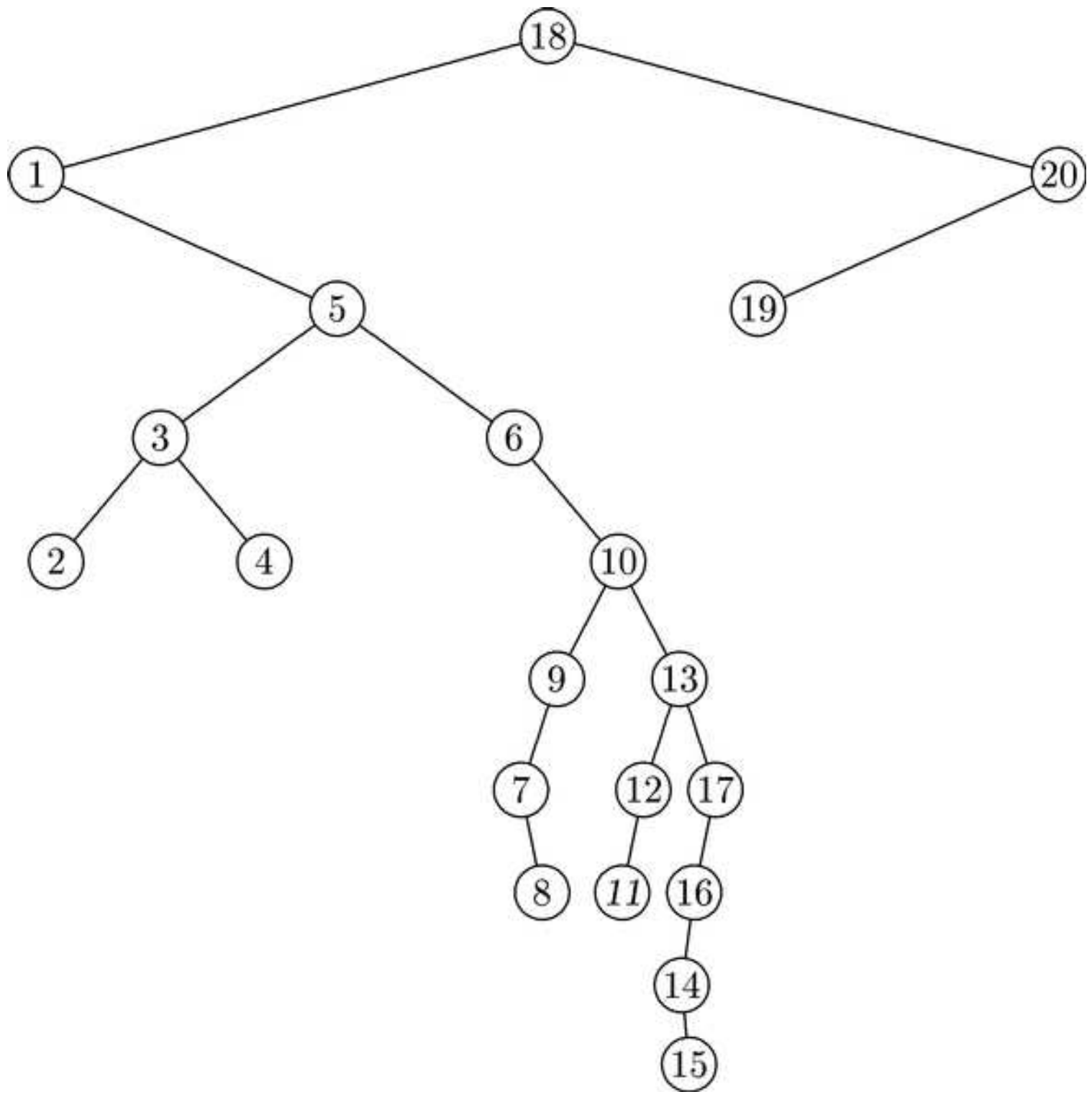

FIG. 1. *The binary search tree with $\pi$ as in Table 1.*

Figure 1 shows the search tree corresponding to the data in Table 1. The crucial point to note is that the path from the root of the tree to the element of interest passes through the descending records in the "+"-list, moving to the left, and the ascending records in the "−"-list, moving to the right.

We recall the definition of some standard distributions: $X$ is said to have a Bernoulli distribution with parameter $p$ if $P(X=1)=1-P(X=0)=p$, to be uniformly distributed on the (finite) set $S$ if $P(X=s)=1/|S|$ for all $s\in S$, and to have a hypergeometric distribution with parameters $N$, $M$

TABLE 1
*A permutation and its subrecord structure for $n=20$, $l=11$*

| | | | | | | | | | | | | | | | | | | | | |
|---|---|---|---|---|---|---|---|---|---|---|---|---|---|---|---|---|---|---|---|---|
| $\pi$ | 18 | 1 | 5 | 6 | 10 | 20 | 3 | 13 | 9 | 17 | 7 | 12 | 8 | 16 | 14 | 19 | 2 | 11 | 4 | 15 |
| $>,<$ | + | − | − | − | − | + | − | + | − | + | − | + | − | + | + | + | − | $\star$ | | |
| records | $\downarrow$ | $\uparrow$ | $\uparrow$ | $\uparrow$ | $\uparrow$ | | | $\downarrow$ | | | | $\downarrow$ | | | | | | $\star$ | | |

and $n$ if

$$P(X=k)=\frac{\binom{M}{k}\binom{N-M}{n-k}}{\binom{N}{n}} \qquad \text{for } k=0,\dots,n.$$

We abbreviate these to $X \sim \mathrm{Ber}(p)$, $X \sim \mathrm{unif}(S)$ and $X \sim \mathrm{HypGeo}(N;M,n)$, respectively. By a random permutation of a finite set $S$ we always mean a permutation that is uniformly distributed on the $|S|!$ possible values.

THEOREM 1. *Suppose that $N, G_{1,l}, \dots, G_{n,l}, K_1, K_2, K_3, \dots, K_1', K_2', K_3', \dots$ are independent random variables with $N \sim \mathrm{unif}(\{1,\dots,n\})$, $G_{m,l} \sim \mathrm{HypGeo}(n-1; l-1, m-1)$ for $m=1,\dots,n$ and $K_i, K_i' \sim \mathrm{Ber}(1/i)$ for all $i \in \mathbb{N}$. Then*

$$X_{n,l} \overset{\mathrm{distr}}{=} \sum_{i=1}^{G_{N,l}} K_i + \sum_{i=1}^{N-1-G_{N,l}} K_i'.$$

PROOF. We first formalize the construction that we outlined above with the help of an example. Remember that $n$ and $l$ are given. Let $\pi$ be a random permutation of $\{1,\dots,n\}$ and let $N := \pi^{-1}(l)$ be the position of $l$. Further, let

$$S_- := \{1 \le i < N : \pi(i) < l\} = \{i_1, \dots, i_G\},$$
$$S_+ := \{1 \le i < N : \pi(i) > l\} = \{j_1, \dots, j_{N-1-G}\},$$

with $i_1 < \cdots < i_G$ and $j_1 < \cdots < j_{N-1-G}$ and let

$$\pi_- := (\pi(i_1), \dots, \pi(i_G)), \qquad \pi_+ := (\pi(j_1), \dots, \pi(j_{N-1-G})),$$
$$R_- := \sum_{r=1}^{G} \prod_{k=1}^{r-1} \mathbb{1}_{\{\pi(i_k) < \pi(i_r)\}}, \qquad R_+ := \sum_{r=1}^{N-1-G} \prod_{k=1}^{r-1} \mathbb{1}_{\{\pi(i_k) > \pi(i_r)\}}.$$

With these constructions we have $X_{n,l} = R_- + R_+$; see Section 13.4 in Cormen, Leiserson and Rivest (1990) for a formal proof. It remains to verify the distributional statements. For these, we simply recall some well-known or easily checked properties of records and random permutations; see, for example, Arnold, Balakrishnan and Nagaraja (1998): Obviously, $N \sim \mathrm{unif}(\{1,\dots,n\})$. Given $N=m$, $(\pi(1),\dots,\pi(m-1))$ is a random permutation of the set $\{\pi(i) : 1 \le i < m\}$. We can view $\pi(i)$ as the result of the $i$th draw, without replacement, from an urn with $n-1$ balls, $l-1$ being "white," meaning a result less than $l$. Hence, conditionally on $N=m$,

$$G := |S_-| \sim \mathrm{HypGeo}(n-1; l-1, m-1).$$

Conditionally on $N=m$ and $G=k$, $\pi_-$ and $\pi_+$ are independent random permutations of $\{\pi(i_1),\dots,\pi(i_k)\}$ and $\{\pi(j_1),\dots,\pi(j_{m-1-k})\}$, respectively. The distributional structure of records in random permutations is such that the products in the definition of $R_-$ and $R_+$, which indicate the presence of a record at position $r$, are independent and Bernoulli distributed with parameter $1/r$.

The assertion of the theorem now follows on comparing the respective (conditional) distributions in the above decomposition to those in the constructive representation. □

From the proof of the theorem it is evident that the first sum in the representation corresponds to the number of moves to the right on the path from the root to the node containing $l$; similarly, the second sum corresponds to the moves to the left. In this context it is interesting to note that

$$G_{N,l} \sim \mathrm{unif}(\{0,\dots,l-1\}), \qquad N-1-G_{N,l} \sim \mathrm{unif}(\{0,\dots,n-l\}).$$

To see this, we simply calculate

$$\begin{aligned} P(G_{N,l}=k) &= \frac{1}{n}\sum_{m=1}^{n} P(G_{m,l}=k) \\ &= \frac{1}{n}\sum_{m=1}^{n} \frac{\binom{l-1}{k}\binom{n-l}{m-1-k}}{\binom{n-1}{m-1}} \\ &= \frac{1}{n}\sum_{m=0}^{n-1} \frac{\binom{m}{k}\binom{n-1-m}{l-1-k}}{\binom{n-1}{l-1}} \\ &= \frac{1}{n}\frac{\binom{n}{l}}{\binom{n-1}{l-1}} = \frac{1}{l} \qquad \text{for } k=0,\dots,l-1, \end{aligned}$$

using $\mathrm{HypGeo}(n-1; l-1, m-1) = \mathrm{HypGeo}(n-1; m-1, l-1)$ and one of the basic identities for binomial coefficients given, for example, as equation (5.26) in Graham, Knuth and Patashnik (1989). The statement on $N-1-G_{N,l}$ follows from similar calculations or from symmetry considerations (see also Section 3).

Note, however, that $G_{N,l}$ and $N-1-G_{N,l}$ are not independent; their joint distribution, which will be used repeatedly below, is given by

$$\text{(JD)} \qquad P(G_{N,l}=i, N-1-G_{N,l}=j) = \frac{1}{n}\frac{\binom{i+j}{i}\binom{n-1-i-j}{l-1-i}}{\binom{n-1}{l-1}}.$$

For our first approximation result we require the following bound for the variance of $H(G)+H(N-1-G)$, where we have written $H(G)$ instead of $H_G$. As usual, we put $H(0)=H_0=0$.

LEMMA 2. *Let $G$ and $N$ be random variables with joint distribution given by* (JD). *Then*

$$\mathrm{var}(H(G)+H(N-1-G)) \le 28.$$

PROOF. Because of $\operatorname{var}(X+Y) \le \operatorname{var}(X)+\operatorname{var}(Y)+2\operatorname{var}(X)^{1/2}\operatorname{var}(Y)^{1/2}$, it is enough to bound the variance of $H(G)$ and $H(N-1-G)$ by 7. The remarks following Theorem 1 imply that both can (individually) be represented in distribution as $H(\lfloor kU \rfloor)$ with $U \sim \operatorname{unif}(0,1)$ and $k = l$ and $k = n-1-l$, respectively. We may assume that $k \ge 1$, and then, using Minkowski's inequality,

$$\begin{aligned}\operatorname{var}(H(\lfloor kU\rfloor)) &\le E(H(\lfloor kU\rfloor) - \log k)^2 \\ &= E((H(\lfloor kU\rfloor) - \log k)\mathbb{1}_{\{U<1/k\}})^2 \\ &\quad + E((H(\lfloor kU\rfloor) - \log k)\mathbb{1}_{\{U\ge 1/k\}})^2 \\ &\le \frac{(\log k)^2}{k} + ((EV_k^2)^{1/2} + (EW_k^2)^{1/2})^2,\end{aligned}$$

with

$$\begin{aligned}V_k &:= (H(\lfloor kU\rfloor) - \log(\lfloor kU\rfloor))\mathbb{1}_{\{U\ge 1/k\}}, \\ W_k &:= \mathbb{1}_{\{U\ge 1/k\}} \log\frac{\lfloor kU\rfloor}{k}.\end{aligned}$$

The first term is bounded by $4e^{-2}$; for $V_k$ we use that $|H_j - \log j| \le 1$ for all $j \in \mathbb{N}$. Finally,

$$EW_k^2 = \frac{1}{k}\sum_{j=1}^{k-1}\left(\log\frac{j}{k}\right)^2 \le \int_0^1 (\log x)^2\,dx = 2,$$

which gives $\operatorname{var}(H(\lfloor kU\rfloor)) \le 4e^{-2} + (1+\sqrt{2})^2 < 7$. $\square$

Our first result shows that the distribution of $X_{n,l}$ can be approximated by a Poisson distribution with the same mean; it comes with an explicit error bound. Recall that the total variation distance of two probability measures $\mu$ and $\nu$ concentrated on $\mathbb{N}_0$ is given by

$$d_{\mathrm{TV}}(\mu,\nu) = \sup_{A\subset\mathbb{N}_0} |\mu(A) - \nu(A)| = \tfrac{1}{2}\sum_{k=0}^{\infty} |\mu(\{k\}) - \nu(\{k\})|.$$

Further, for a probability measure $\nu$ concentrated on the nonnegative half line $[0,\infty)$ we write $\operatorname{MixPo}(\nu)$ for the mixed Poisson distribution with mixing measure $\nu$, that is,

$$\operatorname{MixPo}(\nu)(\{k\}) = \int e^{-\lambda}\frac{\lambda^k}{k!}\nu(d\lambda) \qquad \text{for all } k \in \mathbb{N}_0.$$

With $\nu = \delta_\lambda$, the one-point measure on $\lambda > 0$, we obtain the usual Poisson distribution $\operatorname{Po}(\lambda)$. This also holds for $\lambda = 0$ as we interpret $\operatorname{Po}(0)$ as $\delta_0$.

THEOREM 3. *With the above notation,*

$$\sup_{l\in\{1,\dots,n\}} d_{\mathrm{TV}}(\mathcal{L}(X_{n,l}),\mathrm{Po}(EX_{n,l})) \le \frac{28+\pi^2}{\log n} \qquad \text{for all } n\ge 2.$$

PROOF. We first give a conditional approximation by a Poisson distribution which leads to an approximation by a mixed Poisson distribution. The latter will then be approximated by a Poisson distribution with the same mean.

We use the following fundamental Poisson approximation result: If $X_1,\dots,X_n$ are independent with $X_i \sim \mathrm{Ber}(p_i)$, then

$$d_{\mathrm{TV}}\left(\mathcal{L}\left(\sum_{i=1}^n X_i\right),\mathrm{Po}\left(\sum_{i=1}^n p_i\right)\right) \le \frac{1}{\sum_{i=1}^n p_i}\sum_{i=1}^n p_i^2;$$

see, for example, page 8 in Barbour, Holst and Janson (1992). Together with the representation in Theorem 1 this immediately implies the following bound for the Poisson approximation of the conditional distributions:

$$\begin{aligned} &d_{\mathrm{TV}}(\mathcal{L}(X_{n,l}|G=i,N-1-G=j),\mathrm{Po}(H_i+H_j)) \\ &\qquad \le \frac{1}{H_i+H_j}\left(\sum_{l=1}^i \frac{1}{l^2}+\sum_{l=1}^j\frac{1}{l^2}\right) \\ &\qquad \le \frac{\pi^2}{3H_{i+j}} \end{aligned}$$

for $i+j>0$; for $i=j=0$ the distance is 0. Note that $i+j$ corresponds to $N-1$, which is uniformly distributed on $\{0,\dots,n-1\}$. The unconditioning step therefore leads to

$$d_{\mathrm{TV}}(\mathcal{L}(X_{n,l}),\mathrm{MixPo}(\mu_{n,l})) \le \frac{\pi^2}{3}\frac{1}{n}\sum_{m=1}^{n-1}\frac{1}{H_m},$$

where $\mu_{n,l}$ denotes the distribution of $H(G)+H(N-1-G)$. Standard elementary arguments show that $\sum_{m=1}^{n-1} 1/H_m \le 3n/\log n$ for $n\ge 2$.

A mixed Poisson distribution can be approximated by an ordinary Poisson distribution with the same mean. Using total variation distance we have, according to Theorem 1.C(ii) in Barbour, Holst and Janson (1992),

$$d_{\mathrm{TV}}(\mathrm{MixPo}(\mu_{n,l}),\mathrm{Po}(EX_{n,l})) \le \frac{\sigma^2}{EX_{n,l}},$$

with $\sigma^2$ the variance associated with $\mu_{n,l}$. Here we have used that the expectation associated with $\mu_{n,l}$ is equal to $EX_{n,l}$. An appeal to Lemma 2 and the triangle inequality now completes the proof. □

We do not claim that the numerical values in the bound are tight; for us, the more important aspect is the fact that the bound does not depend on $l$. In particular, with $(l_n)_{n\in\mathbb{N}}$ a sequence of integers with $1 \le l_n \le n$ for all $n \in \mathbb{N}$, but completely arbitrary otherwise, and $Y_n \sim \mathrm{Po}(EX_{n,l_n})$,

$$\left|P\left(\frac{X_{n,l_n} - EX_{n,l_n}}{\sqrt{EX_{n,l_n}}} \le x\right) - \Phi(x)\right|$$
$$\le d_{\mathrm{TV}}(\mathcal{L}(X_{n,l_n}), \mathcal{L}(Y_n)) + \left|P\left(\frac{Y_n - EY_n}{\sqrt{\mathrm{var}(Y_n)}} \le x\right) - \Phi(x)\right|$$

for all $x \in \mathbb{R}$, so the asymptotic normality of Poisson distributions with parameter tending to infinity and the bound in Theorem 3 together imply that

$$\frac{X_{n,l_n} - EX_{n,l_n}}{\sqrt{EX_{n,l_n}}} \stackrel{\mathrm{distr}}{\to} Z, \qquad Z \sim N(0,1),$$

as $n \to \infty$. [In fact, combining this with the Berry–Esseen theorem we obtain the rate $O((\log n)^{-1/2})$ for the Kolmogorov–Smirnov distance.] Special cases can be obtained on using Arora and Dent's formula (AD). For example, if

$$\text{(SC)} \qquad \lim_{n\to\infty} \frac{\min\{\log(l_n), \log(n - l_n)\}}{\log n} = t$$

for some $t \in [0,1]$, then

$$\frac{X_{n,l_n} - (1+t)\log n}{\sqrt{(1+t)\log n}} \stackrel{\mathrm{distr}}{\to} Z, \qquad Z \sim N(0,1).$$

In particular, if $l_n/n \to t \in (0,1)$, then $(X_{n,l_n} - 2\log n)/\sqrt{2\log n}$ is asymptotically standard normal, irrespective of the value of $t$.

Louchard (1987) showed that, with $U_n \sim \mathrm{unif}(\{1,\dots,n\})$ independent of the search trees,

$$\frac{X_{n,U_n} - 2\log n}{\sqrt{2\log n}} \stackrel{\mathrm{distr}}{\to} Z, \qquad Z \sim N(0,1).$$

This can now be derived from (SC) via the representation $U_n = \lceil nU \rceil$ with $U \sim \mathrm{unif}(0,1)$ by conditioning on $U = t \in (0,1)$. (A conditioning argument can also be used to extend the bound in Theorem 3 to randomly chosen $l$-indices.) The special case also makes precise the intuitive picture that nodes with extreme keys, that is, with $l$ being close to 1 or $n$, have lesser depth and will be found faster than those "within" the range from 1 to $n$.

In order to see the influence of the key on the node depth in the midrange, by which we mean that $l_n/n \to t$ for some $t$ with $0 < t < 1$, we have to use a different metric for probability distributions. This becomes obvious as soon we expand $EX_{n,l_n}$ up to constants, since in an asymptotic normality result constant shifts do not matter asymptotically if the scaling factors tend to infinity. If we use the total variation distance, this even holds on the Poisson approximation level as

$$\lim_{\lambda\to\infty} d_{\mathrm{TV}}(\mathrm{Po}(\lambda+c),\mathrm{Po}(\lambda)) = 0 \qquad \text{for all } c>0.$$

Our second result shows that with a suitable Wasserstein metric shifts do become visible. There are two consequences: We now need a mixed Poisson distribution as approximating measure, and we lose on the rate side. Following Barbour, Holst and Janson (1992), we consider the distance $d_W$ for probability distributions $\mu$, $\nu$ on (the Borel subsets of) the real line defined by

$$d_W(\mu,\nu) := \sup\left\{\left|\int f\,d\mu - \int f\,d\nu\right| : f:\mathbb{R}\to\mathbb{R}, \sup_{|x-y|\le 1}|f(x)-f(y)|\le 1\right\}.$$

For distributions concentrated on the nonnegative integers it can be shown that

$$d_W(\mu,\nu) = \sum_{k=0}^{\infty}|\mu([k,\infty)) - \nu([k,\infty))|.$$

Hence, if $X$ and $Y$ are random variables with distributions $\mu$ and $\nu$, respectively, then $d_W(\mu,\nu)\ge |EX-EY|$, which in turn implies that $\mathrm{Po}(\lambda+c)$ and $\mathrm{Po}(\lambda)$ remain distinguishable under this distance if $\lambda\to\infty$, $c>0$ fixed (we generally use $d_W$ only in connection with distributions with finite mean). Further, $d_W$ can be realized by a suitable coupling in the sense that

$$d_W(\mu,\nu) = \min\{E|X-Y| : X\sim\mu, Y\sim\nu\}.$$

The following lemma contains two properties of the Wasserstein distance; their proof makes use of the above alternative expressions for $d_W$. When we use the first of these below we will speak of unconditioning; a similar property for the total variation distance has already been used in the proof of Theorem 3. The second property shows that $\mu\mapsto \mathrm{MixPo}(\mu)$ is a weak $d_W$-contraction.

LEMMA 4. (a) *If $X$ with $P(X\in\mathbb{N})=1$ and $Y$ are random variables such that*

$$d_W(\mathcal{L}(X|Y=y),\mathrm{Po}(\phi(y))) \le f(y)$$

*for all $y$, with measurable functions $\phi$ and $f$, then*

$$d_W(\mathcal{L}(X),\mathrm{MixPo}(\mathcal{L}(\phi(Y)))) \le Ef(Y).$$

(b) *For any two probability distributions $\mu, \nu$ on the nonnegative real line,*

$$d_W(\text{MixPo}(\mu), \text{MixPo}(\nu)) \le d_W(\mu, \nu).$$

PROOF. (a) We condition on the value of $Y$; $\int \cdots \mathcal{L}(Y)(dy)$ means that we integrate with respect to the distribution of $Y$:

$$\begin{aligned}
&d_W(\mathcal{L}(X), \text{MixPo}(\mathcal{L}(\phi(Y))))\\
&\qquad = \sum_{k=0}^{\infty} \left| \int (\mathcal{L}(X|Y=y)([k,\infty)) - \text{Po}(\phi(y))([k,\infty)))\mathcal{L}(Y)(dy) \right|\\
&\qquad \le \int \sum_{k=0}^{\infty} |\mathcal{L}(X|Y=y)([k,\infty)) - \text{Po}(\phi(y))([k,\infty))| \mathcal{L}(Y)(dy)\\
&\qquad = \int d_W(\mathcal{L}(X|Y=y), \text{Po}(\phi(y)))\mathcal{L}(Y)(dy).
\end{aligned}$$

(b) Let $(N_t)_{t\ge 0}$ be a unit rate Poisson process and let $X$ and $Y$ be random variables, independent of the process, with $X \sim \mu$, $Y \sim \nu$ and $d_W(\mu,\nu) = E|X-Y|$. Then $N_X \sim \text{MixPo}(\mu)$, $N_Y \sim \text{MixPo}(\nu)$ so that by conditioning on $X$ and $Y$ and considering the cases $X > Y$ and $X \le Y$ separately,

$$\begin{aligned}
d_W(\text{MixPo}(\mu), \text{MixPo}(\nu)) &\le E|N_X - N_Y|\\
&= E(E[|N_X - N_Y|\,|X, Y])\\
&= E|X-Y| = d_W(\mu,\nu). \qquad \square
\end{aligned}$$

We also need an elementary estimate related to hypergeometric distributions.

LEMMA 5. *With $X \sim \text{HypGeo}(N; M, n)$,*

$$E\left(\left|\log \frac{X}{EX}\right| \mathbb{1}_{\{X>0\}}\right) \le \frac{4N\log N}{nM} + 2\sqrt{\frac{N}{nM}}.$$

PROOF. We use

$$EX = \frac{nM}{N}, \qquad \text{var}(X) \le \frac{nM}{N}$$

together with Chebyshev's inequality, the bound $\log N$ for the integrand, the fact that $|\log(1+x)| \le 2|x|$ on $|x| \le 1/2$, and $E|X - EX| \le \sqrt{\text{var}(X)}$ to obtain

$$\begin{aligned}
E\left(\left|\log \frac{X}{EX}\right| \mathbb{1}_{\{X>0\}}\right) &\le (\log N) P\left(|X - EX| \ge \frac{EX}{2}\right) + 2E\left|\frac{X}{EX} - 1\right|\\
&\le \frac{4N\log N}{nM} + 2\sqrt{\frac{N}{nM}}.
\end{aligned}$$

□

We can now state and prove our second approximation result for key values in the central range.

THEOREM 6. *Suppose that $l_n$ varies with $n$ such that*

$$\frac{l_n}{n} = t + O\left(\frac{1}{\sqrt{\log n}}\right)$$

*with some $t \in (0,1)$. Let*

$$\nu_{n,t} := \mathcal{L}((2\log n + 2\gamma + \log(t(1-t)) - 2X)^+),$$

*where $X$ is exponentially distributed with mean 1. Then*

$$d_W(\mathcal{L}(X_{n,l_n}), \mathrm{MixPo}(\nu_{n,t})) = O\left(\frac{1}{\sqrt{\log n}}\right).$$

PROOF. We continue to use the notation introduced in the proof of Theorem 3 and again begin by comparing conditional distributions to Poisson distributions. The basic result for the Wasserstein distance, obtained by combining Lemma 1.1.5 and Remark 1.1.7 in Barbour, Holst and Janson (1992), is the following: If $X_1, \ldots, X_n$ are independent with $X_i \sim \mathrm{Ber}(p_i)$, then

$$d_W\left(\mathcal{L}\left(\sum_{i=1}^n X_i\right), \mathrm{Po}\left(\sum_{i=1}^n p_i\right)\right) \le \frac{2}{\sqrt{\sum_{i=1}^n p_i}} \sum_{i=1}^n p_i^2.$$

In our situation we obtain with the representation in Theorem 1, abbreviating $G_{N_n,l_n}$ to $G_n$,

$$\begin{aligned} & d_W(\mathcal{L}(X_{n,l_n} | G_n = i, N_n - 1 - G_n = j), \mathrm{Po}(H_i + H_j)) \\ & \quad \le \frac{2}{\sqrt{H_i + H_j}} \left(\sum_{m=1}^{i} \frac{1}{m^2} + \sum_{m=1}^{j} \frac{1}{m^2}\right) \le \frac{2\pi^2}{3\sqrt{H_{i+j}}}. \end{aligned}$$

Unconditioning and using $(\log n)^{1/2} \sum_{m=1}^{n-1} H_m^{-1/2} = O(n)$, we see that

$$d_W(\mathcal{L}(X_{n,l}), \mathrm{MixPo}(\mu_n)) = O\left(\frac{1}{\sqrt{\log n}}\right),$$

where $\mu_n := \mathcal{L}(H(G_n) + H(N_n - 1 - G_n))$. Using the triangle inequality and Lemma 4(b) we see that it remains to show that $d_W(\mu_n, \nu_{n,t}) = O((\log n)^{-1/2})$. This will follow if we can find random variables $X_n$ and $Y_n$ such that $\mathcal{L}(X_n) = \mu_n$, $\mathcal{L}(Y_n^+) = \nu_{n,t}$ and $\sqrt{\log n} E|X_n - Y_n| = O(1)$. (Because of $X_n \ge 0$, going from $Y_n$ to $Y_n^+$ will not increase the Wasserstein distance to the distribution of $X_n$.) Let $U \sim \mathrm{unif}(0,1)$ and $N_n := \lceil nU \rceil$ for all $n \in \mathbb{N}$. With

$$X_n := H(G_n) + H(N_n - 1 - G_n),$$
$$Y_n := 2\log n + 2\gamma + \log(t(1-t)) + 2\log U,$$

the distributional requirements are satisfied and we have

$$|X_n - Y_n| \le \sum_{i=1}^{6} |Z_{i,n}|$$

with

$$Z_{1,n} := H(G_n) - \log(G_n)\mathbb{1}_{\{G_n>0\}} - \gamma,$$
$$Z_{2,n} := \log(G_n)\mathbb{1}_{\{G_n>0\}} - \log(l_n - 1) - \log\Big(\frac{N_n}{n}\Big),$$
$$Z_{3,n} := \log(l_n - 1) - \log(nt) + \log\Big(\frac{N_n}{n}\Big) - \log U,$$
$$Z_{4,n} := H(N_n - 1 - G_n) - \log(N_n - 1 - G_n)\mathbb{1}_{\{N_n-1-G_n>0\}} - \gamma,$$
$$Z_{5,n} := \log(N_n - 1 - G_n)\mathbb{1}_{\{N_n-1-G_n>0\}} - \log(n - 1 - l_n) - \log\Big(\frac{N_n}{n}\Big),$$
$$Z_{6,n} := \log(n - 1 - l_n) - \log(n(1-t)) + \log\Big(\frac{N_n}{n}\Big) - \log U.$$

For the first of these we use the fact that, for some constant $C < \infty$,

$$|H_n - \log n - \gamma| \le \frac{C}{n} \qquad \text{for all } n \in \mathbb{N},$$

and $\mathcal{L}(G_n) = \mathrm{unif}(\{0,\dots,l_n - 1\})$ to obtain

$$E|Z_{1,n}| \le \gamma P(G_n = 0) + CE\Big(\frac{1}{G_n}\mathbb{1}_{\{G_n>0\}}\Big)$$
$$= \frac{\gamma}{l_n} + \frac{C}{l_n}\sum_{k=1}^{l_n-1}\frac{1}{k} = O\Big(\frac{\log n}{n}\Big).$$

The second term is slightly more complicated as it involves both $G_n$ and $N_n$. Conditioning on the latter we get

$$E|Z_{2,n}| \le E(E[|Z_{2,n}||N_n]).$$

On $\{N_n = 1\}$ we have $G_n \equiv 0$, which leads to

$$E[|Z_{2,n}||N_n = 1] = \log\Big(\frac{n}{l_n - 1}\Big) = O(1).$$

Together with $P(N_n = 1) = 1/n$ this gives $E|Z_{2,n}|\mathbb{1}_{\{N_n=1\}} = O(1/n)$. We may therefore assume that $N_n > 1$ as long as we deal with $Z_{2,n}$.

We use another decomposition,

$$|Z_{2,n}| \le Z_{2,1,n} + Z_{2,2,n} + Z_{2,3,n}$$

with

$$\begin{aligned} Z_{2,1,n} &= \left|\log \frac{G_n}{E[G_n|N_n]}\right| \mathbb{1}_{\{G_n>0\}}, \\ Z_{2,2,n} &= \left|\log \frac{N_n - 1}{n-1} - \log \frac{N_n}{n}\right|, \\ Z_{2,3,n} &= \left|\log \frac{N_n(l_n - 1)}{n}\right| \mathbb{1}_{\{G_n=0\}}. \end{aligned}$$

Lemma 5 yields

$$E\left[\left|\log \frac{G_n}{E[G_n|N_n]}\right| \mathbb{1}_{\{G_n>0\}} \Big| N_n\right] \le \frac{4(n-1)\log(n-1)}{(N_n-1)(l_n-1)} + 2\sqrt{\frac{n-1}{(N_n-1)(l_n-1)}}$$

on $N_n > 1$, which together with

$$\frac{1}{n}\sum_{k=1}^{n-1} \frac{\log(n-1)}{k} = O\left(\frac{1}{\sqrt{\log n}}\right), \qquad \frac{1}{n}\sum_{k=1}^{n-1} \frac{1}{\sqrt{k}} = O\left(\frac{1}{\sqrt{\log n}}\right),$$

gives $EZ_{2,1,n}\mathbb{1}_{\{N_n>1\}} = O((\log n)^{-1/2})$. For $Z_{2,2,n}$ we obtain

$$\begin{aligned} EZ_{2,2,n}\mathbb{1}_{\{N_n>1\}} &= \frac{1}{n}\sum_{k=2}^{n} \left|\log \frac{k-1}{n-1} - \log \frac{k}{n}\right| \\ &\le \frac{1}{n}\sum_{k=2}^{n} (\log k - \log(k-1)) + \frac{1}{n}\sum_{k=2}^{n} (\log n - \log(n-1)) \\ &= \frac{1}{n}\log n + \frac{n-1}{n}\log \frac{n}{n-1} \\ &= O\left(\frac{\log n}{n}\right). \end{aligned}$$

On $\{N_n > \sqrt{n}\}$ we have

$$E[\mathbb{1}_{\{G_n=0\}}|N_n] \le \left(\frac{n - l_n}{n - \sqrt{n}}\right)^{\sqrt{n}} \le \kappa^{\sqrt{n}}$$

for some $\kappa < 1$ and $n$ large enough, hence

$$\begin{aligned} EZ_{2,3,n}\mathbb{1}_{\{N_n>1\}} &= E(E[Z_{2,3,n}|N_n]\mathbb{1}_{\{1<N_n\le\sqrt{n}\}}) \\ &\quad + E(E[Z_{2,3,n}|N_n]\mathbb{1}_{\{N_n>\sqrt{n}\}}) \end{aligned}$$

$$\leq \frac{1}{n}\sum_{k=2}^{\lfloor\sqrt{n}\rfloor}\left|\log\frac{k(l_n-1)}{n}\right| + \frac{1}{n}\sum_{k=\lfloor\sqrt{n}\rfloor+1}^{n}\log(l_n)\kappa^{\sqrt{n}}.$$

Both terms on the right-hand side are obviously $O((\log n)^{-1/2})$ so that this rate also holds for $EZ_{2,3,n}\mathbb{1}_{\{N_n>1\}}$ and therefore for $E|Z_{2,n}|$ too.

For $Z_{3,n}$ we use the rate condition on $\frac{l_n}{n} - t$ together with the following argument which is based on the construction of $N_n$:

$$\begin{aligned} E\left|\log\frac{N_n}{n} - \log U\right| &= E\log\frac{N_n}{n} - E\log U \\ &= \frac{1}{n}\sum_{k=1}^{n}\log\frac{k}{n} + 1 \\ &= \frac{1}{n}\log(n!) - \log n + 1 \\ &= O\left(\frac{\log n}{n}\right). \end{aligned}$$

Finally, adapting the arguments used for $Z_{i,n}$ to $Z_{i+3,n}$, $i = 1, 2, 3$, is a straightforward task. □

**3. Miscellaneous comments.** We relate our findings to another classical algorithm in Section 3.1. In Section 3.2 we discuss the expectation and the variance of $X_{n,l}$. The use of (and need for) other probability metrics, together with the relationship between total variation and Wasserstein distance, are briefly considered in Section 3.3. The final subsection deals with another noteworthy aspect of the representation of $X_{n,l}$ as the sum of the number of moves to the right and the number of moves to the left.

3.1. A situation very similar to the one considered above arises in connection with Hoare's (1961) selection algorithm FIND, a randomized divide-and-conquer algorithm that selects the $l$th smallest element of a totally ordered set $S$ of size $n$ in a recursive manner: First, an $x$ from $S$ is chosen uniformly at random. Comparing this element to all others, we obtain the subsets $S_- := \{y \in S : y < x\}$ and $S_+ := \{y \in S : y > x\}$. We continue with $(l, S)$ replaced by $(l, S_-)$ if the size $k := |S_-|$ is greater than or equal to $l$ and with $(l-1-k, S_+)$ if $k < l-1$. If $k = l-1$, then we stop and return $x$. For the time required by the algorithm the number of comparisons $C_{n,l}$ is most important, but the number $R_{n,l}$ of recursions has also been investigated. Instead of introducing randomness via the selection of the pivotal element, we can equivalently assume that the data are random, with all permutations being equally likely, that we operate on lists rather than sets and that we always choose the first element of the list as the pivot. This connects

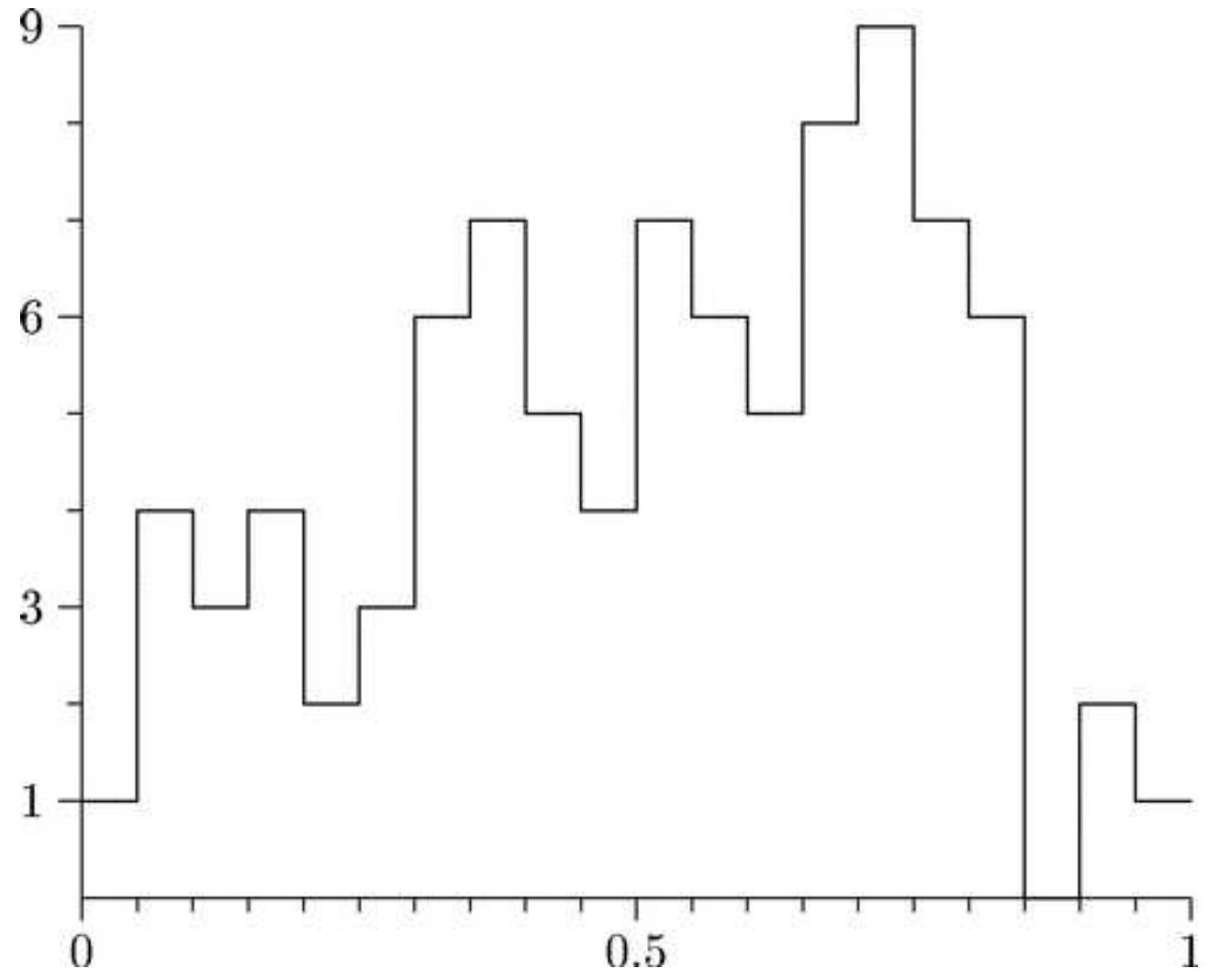

Fig. 2. *The depth plot* $t \mapsto X_{n,\lceil nt \rceil}$ *for the permutation in Table* 1.

Find to binary search trees, with $S_-$ and $S_+$ corresponding to the left and right subtree, respectively, and indeed, it is well known that $R_{n,l}$ is equal in distribution to $X_{n,l}$ (or to $1 + X_{n,l}$ if we include the initial step).

Again, details are given in many of the standard textbooks; see also the recent book by Mahmoud (2000). As with binary search trees, if interest is in the behavior of these quantities for $n$ large, one can average out the $l$. This leads to results on the number of comparisons and recursions needed for a randomly chosen $l$; see, for example, Section 7.5 in Mahmoud (2000) and the references given there. Instead, Grübel and Rösler (1996) considered the whole function $l \mapsto C_{n,l}$. The resulting limit theorem for the stochastic processes $(C_{n,\lceil tn \rceil})_{0 \le t \le 1}$ implies the distributional convergence of $C_{n,l_n}/n$ if the sequence $(l_n)_{n \in \mathbb{N}}$ is such that $l_n/n \to t$ as $n \to \infty$ for some $t \in [0,1]$; the limit distribution depends on $t$. A different approach, leading to this result more easily, is given in Grübel (1998). The results in the previous section cover similar aspects for the number of recursions required. In particular, the terms $\sum_{i=1}^{G_{N,l}} K_i$ and $\sum_{i=1}^{N-1-G_{N,l}} K_i'$ in Theorem 1 represent the number of times that the element of interest is put into $S_-$ and $S_+$, respectively, in the course of the algorithm. It is interesting to note that, in contrast to the situation with the number of comparisons, we have concentration of mass for the number of recursions in the sense that $R_{n,l_n}/ER_{n,l_n}$ converges to 1 in probability. An analogue to the result in Grübel and Rösler (1996) would be a functional limit theorem for the "depth plot" $l \mapsto X_{n,l}$ which, incidentally, characterizes the binary search tree. Figure 2 shows this plot for the permutation in Table 1.

3.2. The representation in Theorem 1 leads to an alternative proof of (AD). Let $Y := \sum_{i=1}^{G_{N,l}} K_i$, $Z := \sum_{i=1}^{N-1-G_{N,l}} K_i'$, with the notation as in Theorem 1. Using $G_{N,l} \sim \mathrm{unif}(\{0,\dots,l-1\})$, $EK_i = 1/i$ and equation (6.67) in Graham, Knuth and Patashnik (1989), we obtain

$$EY = \frac{1}{l}\sum_{j=0}^{l-1}\sum_{i=1}^{j}\frac{1}{i} = \frac{1}{l}\sum_{j=0}^{l-1} H_j = H_l - 1.$$

Together with a similar calculation for $Z$, this gives

$$EX_{n,l} = EY + EZ = H_l + H_{n+1-l} - 2.$$

The variance of $X_{n,l}$ is mentioned in Arora and Dent (1969); the explicit formula

$$\text{(KP)}\qquad \begin{aligned} \mathrm{var}(X_{n,l}) = {} & \frac{2(n+1)}{l(n+1-l)} H_n + \left(1 - \frac{2(n+1)}{l(n+1-l)}\right)(H_l + H_{n+1-l}) \\ & - H_l^{(2)} - H_{n+1-l}^{(2)} + \frac{2}{l(n+1-l)} + 2, \end{aligned}$$

with $H_n^{(2)} := \sum_{k=1}^{n} 1/k^2$, is given in Kirschenhofer and Prodinger (1998). Obtaining this from our representation is a somewhat tedious task that boils down to an unsightly formula involving harmonic numbers and a multitude of binomial coefficients. In contrast to the situation with $EX_{n,l}$, this does not seem to lead to an intuitive or short proof, so we do not give the details.

3.3. We have pointed out in Section 2 that the total variation distance will not distinguish between, say, $\mathrm{Po}(\lambda)$ and $\mathrm{Po}(\lambda + c)$ with $c$ constant as $\lambda \to \infty$, so we may have $d_{\mathrm{TV}}(\mathcal{L}(X_n), \mathcal{L}(Y_n)) \to 0$ even if $EX_n - EY_n$ does not vanish asymptotically as $n \to \infty$. For general distributions on the real line we may conversely have a small Wasserstein distance together with a large total variation distance, but for distributions concentrated on the integers the simple relation

$$\mu(\{k\}) = \mu([k,\infty)) - \mu([k+1,\infty))$$

implies that

$$d_{\mathrm{TV}}(\mu,\nu) \le 2 d_W(\mu,\nu).$$

Using $d_W$ instead of $d_{\mathrm{TV}}$, we obtained an approximation that is asymptotically correct with respect to first moments in the sense that $\lim_{n\to\infty} d_W(\mathcal{L}(X_n), \mathcal{L}(Y_n)) = 0$ implies $\lim_{n\to\infty}(EX_n - EY_n) = 0$. From (KP) and some straightforward calculations it follows that we would need yet another metric and a more detailed expansion to obtain an approximation that is asymptotically correct for second moments too; see, for example, the metric used in Mahmoud and Neininger (2003).

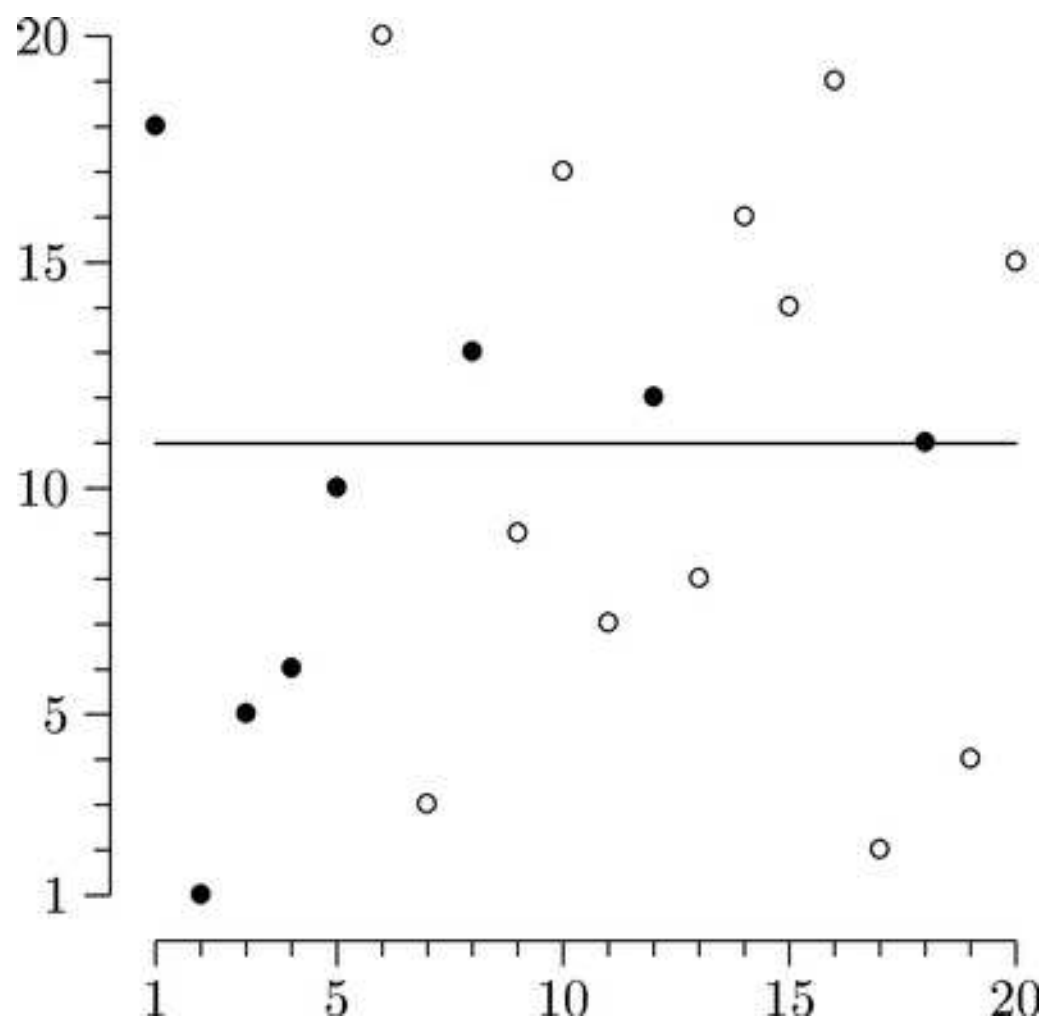

FIG. 3. *The scatterplot for the permutation in Table 1 (•: records in the subpermutations).*

3.4. The simplification for the two constituent parts of $X_{n,l}$ that we used in Section 3.2 has the following noteworthy consequence: The distribution of $\sum_{i=1}^{G_{N,l}} K_i$, with the assumptions as in Theorem 1, is equal to that of $\sum_{i=1}^{l} K_i - 1$, which makes the random summation index disappear. With $X_{n,l}^{\leftarrow}$ and $X_{n,l}^{\rightarrow}$ for the number of moves to the left and right, respectively, this means that

$$\mathcal{L}(1+X_{n,l}^{\rightarrow}) = \mathcal{L}\left(\sum_{i=1}^{l} K_i\right), \qquad \mathcal{L}(1+X_{n,l}^{\leftarrow}) = \mathcal{L}\left(\sum_{i=1}^{n+1-l} K_i\right),$$

with $K_1, K_2, \ldots$ independent and $K_i \sim \mathrm{Ber}(1/i)$. Since $X_{n,l} = X_{n,l}^{\leftarrow} + X_{n,l}^{\rightarrow}$, this leads to another proof of (AD).

A glance at Figure 3 explains the "distributional coincidence": $1 + X_{n,l}^{\rightarrow}$ is the number of ascending records in the subpermutation of the $l$ elements that are less than or equal to $l$, $1 + X_{n,l}^{\leftarrow}$ is the number of descending records in the subpermutation of the $n+1-l$ elements that are greater than or equal to $l$. This leads to a very simple description of the node depth distribution in the extreme cases,

$$\mathcal{L}(1+X_{n,1}) = \mathcal{L}(1+X_{n,n}) = \mathcal{L}\left(\sum_{i=1}^{n} K_i\right),$$

since for the minimum and maximum all steps are in one direction only. Note, however, that despite the independence of the subpermutations of the elements that are *strictly* smaller, respectively larger, than $l$, $X_{n,l}^{\leftarrow}$ and $X_{n,l}^{\rightarrow}$ are *not* independent. Indeed, since $\mathcal{L}(X^{\rightarrow}_{n,l}) = \mathcal{L}(X_{l,l})$ and $\mathcal{L}(X^{\leftarrow}_{n,l}) = \mathcal{L}(X_{n+1-l,1})$, it is tempting to think of $X_{n,l}$ as the sum of $X_{l,l}$ and $X_{n+1-l,1}$, but the simplest nontrivial case already provides a counterexample to the assumption that these can be taken to be independent: $\mathcal{L}(X_{3,2}) = \text{unif}(\{0,1,2\})$, $\mathcal{L}(X_{2,1}) = \mathcal{L}(X_{2,2}) = \text{unif}(\{0,1\})$.

**Acknowledgment.** We thank the referee for drawing our attention to Devroye and Neininger (2004).

Institut für Mathematische Stochastik
Universität Hannover
Postfach 60 09
D-30060 Hannover
Germany
e-mail: rgrubel@stochastik.uni-hannover.de